\documentclass[11pt,a4paper]{article}
\usepackage[sc,noBBpl]{mathpazo}
\usepackage[mathscr]{eucal}
\usepackage{amsmath,amsfonts,amssymb,amsthm}
\newcommand\scG{{\mathscr G}}

\newcommand\mvector{\boldsymbol}

\newcommand\vd{\mvector{d}}

\newcommand\vp{\mvector{p}}
\newcommand\vq{\mvector{q}}

\newcommand\vx{\mvector{x}}
\newcommand\vy{\mvector{y}}
\newcommand\vz{\mvector{z}}

\newcommand\vC{\mvector{C}}

\newcommand\vLambda{\mvector{\Lambda}}

\newcommand\vGamma{\mvector{\Gamma}}

\newcommand\field{\mathbb}

\newcommand\C{\field{C}}
\newcommand\Z{\field{Z}}
\newcommand\M{\field{M}}

\newcommand\bbP{\mathbb{P}}

\newcommand\tr{\operatorname{Tr}}

\newcommand\ord{\operatorname{ord}}
\newcommand\grad{\operatorname{grad}}

\newcommand\trdeg{\operatorname{tr.\!deg}}

\newcommand\rmd{\mathrm{d}}

\newcommand\CP{\ensuremath{\C\bbP}}
\newcommand\rmi{\mathrm{i}\mspace{1mu}}
\newcommand\rme{\mathrm{e}}
\newcommand\Dt{\frac{\mathrm{d}\phantom{t} }{\mathrm{d}\mspace{1mu} t}}

\newcommand\Dtt{\frac{\mathrm{d}^2\phantom{t} }{\mathrm{d}t^2}}

\newcommand\mtext[1]{\quad\text{#1}\quad}

\newcommand\defset[2]{\left\{{#1}\;\vert \;\; {#2} \,\right\}}

\usepackage[square]{natbib}
\theoremstyle{plain}

\newtheorem{theorem}{Theorem}
\newtheorem{lemma}[theorem]{Lemma}
\newtheorem{proposition}[theorem]{Proposition}
\newtheorem{corollary}[theorem]{Corollary}
\theoremstyle{definition}

\newtheoremstyle{note}{\topsep}{\topsep}{\slshape}{}{\scshape}{}{ }{}
\theoremstyle{note}

\numberwithin{equation}{section}
\numberwithin{theorem}{section}

\begin{document}
\thispagestyle{empty}
\vspace*{1em}
\begin{center}
\LARGE\textbf{Integrability of  natural Hamiltonian systems with homogeneous potentials of degree zero}
\end{center}
\vspace*{0.5em}
\begin{center}
\large Guy Casale$^1$, Guillaume Duval$^2$, Andrzej J.~Maciejewski$^3$
 and Maria Przybylska$^{4}$
\end{center}
\vspace{2em}
\hspace*{2em}
\begin{minipage}{0.8\textwidth}
\small
$^1$ IRMAR UMR 6625, Universit\'e de Rennes 1, Campus de Beaulieu,\\
35042 Rennes Cedex, France
(e-mail: guy.casale@univ-rennes1.fr)\\[0.5em]
$^2$1 Chemin du Chateau , 76 430 Les Trois Pierres, France. \\
(e-mail: dduuvvaall@wanadoo.fr)\\[0.5em]
$^3$ Institute of Astronomy,
  University of Zielona G\'ora,
  Licealna 9, \\\quad PL-65--417 Zielona G\'ora, Poland
  (e-mail: maciejka@astro.ia.uz.zgora.pl)\\[0.5em]
$^4$Toru\'n Centre for Astronomy,
  N.~Copernicus University,
  Gagarina 11,\\\quad PL-87--100 Toru\'n, Poland,
  (e-mail: Maria.Przybylska@astri.uni.torun.pl)
\end{minipage}\\[1.5em]
\begin{center}
\today
\end{center}
\begin{abstract}
 We derive necessary conditions for integrability in the Liouville sense of natural Hamiltonian systems with homogeneous potential of degree zero. We derive these conditions through an analysis of the differential Galois group of variational equations along a particular solution generated by a non-zero solution $\vd\in\C^n$ of  nonlinear equations $\grad V(\vd)=\vd$. We proved that if the system integrable then the Hessian matrix $V''(\vd)$ has only integer eigenvalues and is semi-simple.
\end{abstract}

\section{Introduction}
\label{sec:int}
In this paper we consider natural Hamiltonian systems with $n$ degrees of
 freedom for which the Hamiltonian function is of the form
\begin{equation}
 \label{eq:ham}
   H=\frac{1}{2}\sum_{i=1}^np_i^2 +V(\vq),
\end{equation}
where $\vq:=(q_1,\ldots,q_n)$, $\vp:=(p_1,\ldots,p_n)$ are the canonical
coordinates, and $V$ is a homogeneous function of degree $k\in\Z$. Although
systems of this form coming from physics and applied sciences are real we
consider them as complex ones, i.e., we assume that they are defined on complex
symplectic manifold $M=\C^{2n}$ equipped with the canonical symplectic form
\begin{equation*}
 \omega=\sum_{i=1}^n \rmd q_i \wedge \rmd p_i.
\end{equation*}
Thus, Hamilton's equations have the canonical form
\begin{equation}
\label{heq}
 \Dt \vq = \vp, \qquad \Dt  \vp = - V'(\vq),
\end{equation}
where $V'$ denotes the gradient of $V$.
Moreover,  in our setting the time variable $t$ is complex.

Let us assume that the following equation
\begin{equation}
 V'(\vq)=\vq,\mtext{where} V'(\vq):=\grad V(\vq),
\end{equation}
has a non-zero solution $\vd\in\C^n$.  It is called a proper Darboux point of the potential,  and  it defines  a two dimensional plane
\begin{equation}
\label{Pi}
 \Pi(\vd):=\defset{(\vq,\vp)\in\C^{2n}}{(\vq,\vp)=(q \vd,p\vd),\ (q,p)\in\C^2},
\end{equation}
 which is   invariant with respect to  system~\eqref{heq}.  Equations~\eqref{heq} restricted to $\Pi(\vd)$ have the form
\begin{equation}
\label{heq1}
 \dot q =p, \quad \dot p = -q^{k-1}.
\end{equation}
For $k\in\Z^\star$, the phase curves of this one degree of freedom Hamiltonian system are
\begin{equation}
 \Gamma_{k,\varepsilon} :=\defset{(q,p)\in\C^2}{ \frac{1}{2}p^2 +\frac{1}{k}q^k=\varepsilon}\subset\C^2, \qquad \varepsilon \in\C.
\end{equation}
Thus, a solution $(q,p)=(q(t), p(t))$ of \eqref{heq1} gives rise a solution $(\vq(t),\vp(t)):=(q\vd,p\vd)$ of equations~\eqref{heq} with the corresponding phase curve
\begin{equation}
 \vGamma_{k,\varepsilon} :=\defset{(\vq,\vp)\in\C^{2n}}{ (\vq,\vp)=(q \vd,p\vd),\ (q,p)\in  \Gamma_{k,\varepsilon} }\subset\Pi(\vd).
\end{equation}

In~\cite{Morales:01::a} J.~J.~Morales-Ruiz and J.~P.~Ramis analysed the
integrability of described Hamiltonian systems.  To this end they investigated the variational equations which are the linearization of system~\eqref{heq} along  phase curve  $\vGamma_{k,\varepsilon}$ with $\varepsilon \neq 0$. They proved the following theorem.
\begin{theorem}[Morales-Ramis]
\label{thm:mo2}
Assume that the Hamiltonian system defined by  Hamiltonian \eqref{eq:ham} with a
homogeneous
potential  $V\in\C(\vq)$ of degree $k\in\Z^\star$ satisfies the following
conditions:
\begin{enumerate}
 \item there exists a non-zero $\vd\in\C^n$ such that $V'(\vd)=\vd$,
\item the system is integrable in the Liouville sense with first integrals which
are meromorphic in a connected neighbourhood of  $\vGamma_{k,\varepsilon}$, with
$\varepsilon\in\C^\star$.
\end{enumerate}
Then for each eigenvalue $\lambda$ of Hessian matrix $V''(\vd)$, pair
$(k,\lambda)$ belongs to an item of the following list
\begin{equation}
 \label{eq:tab}
\begin{array}{crcr}
\text{case} &  k & \lambda  &\\[0.5em]
\hline
\vbox to 1.3em{}1. &  \pm 2 &  \text{arbitrary} & \\[0.7em]
2.  & k & p + \dfrac{k}{2}p(p-1) & \\[0.9em]
3. & k & \dfrac 1 {2}\left(\dfrac {k-1} {k}+p(p+1)k\right)  & \\[0.9em]
4. & 3 &  -\dfrac 1 {24}+\dfrac 1 {6}\left( 1 +3p\right)^2, & -\dfrac 1
{24}+\dfrac 3 {32}\left(  1  +4p\right)^2 \\[0.7em]
 & & -\dfrac 1 {24}+\dfrac 3 {50}\left(  1  +5p\right)^2,  &
-\dfrac1{24}+\dfrac{3}{50}\left(2 +5p\right)^2 \\[0.9em]
5. & 4 & -\dfrac 1 8 +\dfrac{2}{9} \left( 1+ 3p\right)^2 & \\[0.9em]
6. & 5 & -\dfrac 9 {40}+\dfrac 5 {18}\left(1+ 3p\right)^2, & -\dfrac 9
{40}+\dfrac 1 {10}\left(2
  +5p\right)^2\\[0.9em]
7. & -3 &\dfrac {25} {24}-\dfrac 1 {6}\left( 1 +3p\right)^2, & \dfrac {25}
{24}-\dfrac 3 {32}\left(1 +4p\right)^2 \\[0.9em]
 & & \dfrac {25} {24}-\dfrac 3 {50}\left(1+ 5p\right)^2, & \dfrac {25}
{24}-\dfrac{3}{50}\left(2+ 5p\right)^2\\[0.9em]
8. & -4 & \dfrac 9 8-\dfrac{2}{9}\left( 1+ 3p\right)^2 & \\[0.9em]
9. & -5 & \dfrac {49} {40}-\dfrac {5} {18}\left(1+3p\right)^2 ,& \dfrac {49}
{40}-\dfrac {1} {10}(2 +5p)^2\\[0.9em]
\hline
\end{array}
\end{equation}
where $p$ is an integer.
\end{theorem}
The above theorem is one of the most beautiful applications of the differential
Galois approach to the integrability studies---the so called Morales-Ramis
theory, see
\cite{Audin:01::c,Audin:02::c,Churchill:95::a,Churchill:96::b,Morales:99::c,Morales:01::b1,%
Morales:01::b2}.
The aim of this paper is to use this framework in order to find necessary
conditions for the integrability of homogeneous potentials with degree of
homogeneity $k=0$ which are excluded by assumptions of Theorem~\ref{thm:mo2}.

Our main result is the following theorem.
\begin{theorem}
\label{thm:we}
Assume that the Hamiltonian system defined by  Hamiltonian \eqref{eq:ham} with a
homogeneous
potential  $V\in\C(\vq)$ of degree $k=0$ satisfies the following
conditions:
\begin{enumerate}
 \item there exists a non-zero $\vd\in\C^n$ such that $V'(\vd)=\vd$, and
\item the system is integrable in the Liouville sense with rational  first integrals.
\end{enumerate}
Then:
\begin{enumerate}
\item all eigenvalues of $V''(\vd)$ are integers, and
 \item matrix $V''(\vd)$ is semi-simple.
\end{enumerate}
\end{theorem}
The fact that new obstructions for the integrability appear when the Hessian
matrix $V''(\vd)$ is not semi-simple was observed  recently in
\cite{Maciejewski:09::} where the following theorem was proved.
\begin{theorem}[Duval, Maciejewski]
  Let $V (\vq)$ be a homogeneous potential degree $k \in \Z \setminus \{- 2, 0,
  2\}$, such that there exists a non-zero solution $\vd\in\C^n$ of equation
  $V'(\vd)=\vd$. If the Hamiltonian system generated by~\eqref{eq:ham} is
  integrable in the Liouville sense with first integrals which are meromorphic
  in a connected neighbourhood of $\vGamma_{k,\varepsilon}$, with $\varepsilon\in\C^\star$, then
  the Hessian matrix $V'' (\vd)$ satisfies the following conditions:
  \begin{enumerate}
    \item For each eigenvalue $\lambda$ of $V'' (\vd)$, the pair $(k, \lambda)$
    belongs to Table~\eqref{eq:tab}.
    \item The matrix $V'' (\vd)$ does not have any Jordan block of size $m \geq 3$.
    \item If $V'' (\vd)$ has a Jordan block of size $m = 2$  with corresponding
      eigenvalue $\lambda$, then   $(k, \lambda)$ belongs to a row in
      Table~\eqref{eq:tab} with number greater than two.
  \end{enumerate}
\end{theorem}
In the above statement,  by a  Jordan block of size $m$ with the  eigenvalue
$\lambda$, we mean that the Jordan form  of  $V''(\vd)$ contains a block of the
form
\begin{equation}
B(\lambda, m):= \begin{bmatrix}
                \lambda & 0 &0& \hdotsfor{2}& 0 \\
                 1& \lambda & 0& \hdotsfor{2}&0\\
                 0& 1&\lambda& \hdotsfor{2}&0 \\
                    \hdotsfor{6}\\
                 0& 0& 0& \ldots &1& \lambda
               \end{bmatrix}
\in\M(m,\C),
\end{equation}
where $\M(m,\C)$ denotes the set of $m\times m$ complex matrices.

In the next Section we derive the variational equations  The proof of Theorem~\ref{thm:we} is contained in Section~\ref{sec:r2} and Section~\ref{sec:r2}, where we investigate the differential Galois groups of subsystems of the variational equations. The last Section contains an example of application of Theorem~\ref{thm:we} to two dimensional potentials.

\section{Variational equations}

Let us assume that a non-zero $\vd\in\C^n$ satisfies $V'(\vd)=\vd$, and $k=0$. Then, $\Pi(\vd)$ defined by~\eqref{Pi} is invariant with respect to the system~\eqref{heq}. However,  for $k=0$, the phase curve corresponding to a solution $(q,p)=(q(t),p(t))$ of   one degree of freedom Hamiltonian system~\eqref{heq1} is not algebraic. In fact, for $k=0$, the phase curves of ~\eqref{heq1} are given by.
\begin{equation}
\label{eq:cur}
 \Gamma_{\varepsilon} :=\defset{(q,p)\in\C^2}{ \frac{1}{2}p^2 +\ln q =\varepsilon}, \quad \varepsilon\in\C.
\end{equation}
A particular solution $(q,p)=(q(t),p(t))$ of~\eqref{heq1} which lies on $\Gamma_{\varepsilon}$ gives
a particular solution $(\vq(t),\vp(t)):=(q\vd,p\vd)\in\Pi(\vd)$ of~\eqref{heq} which lies on phase curve
\begin{equation}
 \vGamma_{\varepsilon} :=\defset{(\vq,\vp)\in\C^{2n}}{ (\vq,\vp)=(q \vd,p\vd),\ (q,p)\in  \Gamma_{\varepsilon} }\subset\Pi(\vd).
\end{equation}

The variational equations along $\vGamma_\varepsilon$ have the form
\begin{equation*}
 \Dt \vx  = \vy , \qquad \Dt \vy = -\frac{1}{q^{2}}V''(\vd)\vx,
\end{equation*}
or simply
\begin{equation}
  \Dtt \vx = -\frac{1}{q^{2}}V''(\vd)\vx.
\end{equation}
We can make a linear transformation $\vx =\vC\vz$ such that   the transformed equations take the form
\begin{equation}
  \Dtt \vz = -{q^{k-2}}\vLambda\vz,
\end{equation}
and matrix $\vLambda=\vC^{-1}V''(\vd)\vC$ is the Jordan form of $V''(\vd)$.

From now on we work with a fixed value $\varepsilon\in\C$. In order to apply the differential Galois theory we have to
introduce an appropriate differential field of functions defined on $
\Gamma_{\varepsilon}$. We assumed that the considered first integrals are
rational functions, i.e., elements of the field $\C(\vq,\vp)$. This is why we
choose as our base field the restriction of field $\C(\vq,\vp)$ to $
\vGamma_{\varepsilon}$. The restriction of $\C(\vq,\vp)$ to $\Pi(\vd)$ gives
field $\C(q,p)$, which together with derivation $\Dt$ defined by
\begin{equation}
\label{eq:eqr}
 \Dt q= p, \quad \Dt p=-\frac{1}{q},
\end{equation}
is a differential field. The restriction of $\C(q,p)$ to $ \Gamma_{\varepsilon}$
gives the field
\begin{equation*}
 \C( \exp(\varepsilon -p^2/2),p)\simeq \C(p,\exp(p^2/2)).
\end{equation*}
This field equipped with derivation $\dot p = -\rme^{-\varepsilon}\exp(p^2/2)$
is our base differential field. From now on in all places $q$ denotes
$\exp(\varepsilon -p^2/2)$.

Recall that according the main theorem of the Morales-Ramis theory, if the
system is integrable in the Liouville sense, then the differential Galois group
of the variational equations along a particular solution is virtually Abelian,
i.e., the identity component of this group is Abelian. Thus, we have to check
whether the differential Galois group $\scG$ over the field $\C(q,p)$ of system
\begin{equation}
  \Dtt \vz = -\frac{1}{q^2}\vLambda\vz,
\end{equation}
is virtually Abelian.  Notice that for each eigenvalue $\lambda$ of $V''(\vd)$,
the above system contains as a subsystem equation of the form
\begin{equation}
\label{eq:s1}
 \Dtt x = -\frac{\lambda}{q^2} x,
\end{equation}
and, if $\vLambda$ has a Jordan block with the corresponding eigenvalue
$\lambda$, then it contains the following subsystem
\begin{align}
\label{s12}
 \Dtt x &= -\frac{\lambda}{q^2} x, \\
\label{s13}
\Dtt y &= -\frac{\lambda}{q^2} y - \frac{1}{q^2} x.
\end{align}
To make the above observations useful we invoke the following fact.  The
differential Galois group of a subsystem  is the quotient of the differential
Galois group of the system. This implies that if  the differential Galois group
of a system  is virtually Abelian, then   the differential Galois group of  its
subsystem  is  virtually Abelian, see Section~1.4 in \cite{Maciejewski:09::}.
Thus, to find obstructions to the integrability of the considered systems it is
enough to investigate the above subsystems.

It is convenient to give another form of \eqref{eq:s1} and  \eqref{s12}.  We
take $p$ as  independent variable in those systems.  We have
\begin{equation*}
 \Dt = \dot p \frac{\rmd \phantom{p}}{\rmd p} \mtext{and} \Dtt =
 \ddot p \frac{\rmd \phantom{p}}{\rmd p} +{ \dot p}^2  \frac{\rmd^2 \phantom{p}}{\rmd p^2}.
\end{equation*}
Thus, using \eqref{eq:eqr} we obtain
\begin{equation*}
 \Dtt = \frac{p}{q^2} \frac{\rmd \phantom{p}}{\rmd p} +
 \frac{1}{q^2} \frac{\rmd^2 \phantom{p}}{\rmd p^2}.
\end{equation*}
So, equation~\eqref{eq:s1} transforms into
\begin{equation}
\label{eq:p1}
  x'' +p  x' + \lambda x = 0,
\end{equation}
where prime denotes derivation with respect to $p$.
System~\eqref{s12}-\eqref{s13} transforms into
\begin{equation}
\label{eq:p12}
\left.
\begin{split}
 &  x'' +p  x' + \lambda x = 0, \\
 &  y'' +p  y' + \lambda y +x = 0.
\end{split}
\quad\right \}
\end{equation}
We underline that  considering both of the above systems we have to determine
their differential Galois groups over the field $\C(q,p)$, not over $\C(p)$!

\section{Rank 2 subsystems}
\label{sec:r2}
Let $G$ denotes the differential Galois group of equation~\eqref{eq:p1} over
$(\C(p), \frac{\rmd \phantom{p}}{\rmd p})$. We show the following.
\begin{proposition}
\label{pro:rd}
 If $\lambda \not\in \Z$, then $G=\mathrm{GL}(2,\C)$.  If $\lambda \in \Z$, then
 $G= \C^\ast \ltimes\C$.
\end{proposition}
\begin{proof}
Let $z=-{p^2}/{2}$. Then, after this change of independent variable,
equation~\eqref{eq:p1} is of the form
\begin{equation}
 z  \frac{\rmd^2 x}{\rmd z^2} + (c-z) \frac{\rmd x}{\rmd z} +a x =
 0  \mtext{where} c = \frac{1}{2}, \quad a = -\frac{\lambda}{2}.
\end{equation}
This is the confluent hypergeometric equation in the Kummer
form~\cite{Erdelyi:81::a}. Its differential Galois group over $\C(z)$ was
investigated in~\cite{MR1038060}, see also~\cite{MR1115231,MR1325757}.  From those
investigations we know that
\begin{itemize}
\item if $\lambda \not\in \Z$, then the Galois group is $\mathrm{GL}(2,\C)$,
\item if $\lambda \in \Z$, then Galois group is $\C^\ast \ltimes\C$.
\end{itemize}

This result gives the Galois group of equation~\eqref{eq:p1} over the field
$\C(p^2)$, hence over $\C(p)$.
\end{proof}
For an integer $\lambda$ we can characterise solutions of equation~\eqref{eq:p1}
in the following proposition.
\begin{proposition}
\label{pro:s}
 Let $\lambda\in\Z$. Then
\begin{equation}
 x_\lambda:=\begin{cases}
                       q H_{\lambda-1}(p), &\text{for}\ \lambda\geq 1, \\
                        H_{-\lambda}(-\rmi p) &\text{for}\ \lambda\leq 0,
\end{cases}
\end{equation}
where $H_n$ is the Hermite polynomial of degree $n$, is a solution of
equation~\eqref{eq:p1}. Its other solution is given by
\begin{equation*}
\widetilde x:= x_\lambda \int\frac{q}{x^2_\lambda}\rmd p.
\end{equation*}
\end{proposition}
\begin{proof}
If we put $x=qy$, then equation~\eqref{eq:p1} transforms into
\begin{equation*}
 y'' -py' + n y = 0 \mtext{where} n=\lambda-1.
\end{equation*}
This is the Hermite differential equation which has for integer $n=\lambda-1$
polynomial solution $H_n(p)$.  If $\lambda\leq 0$  we make transformation
$p\mapsto -\rmi p$, which transforms  equation~\eqref{eq:p1} into the Hermite
equation with $n=-\lambda$. The formula for the second solution is standard.
\end{proof}
We notice here that the Hermite polynomials $H_n$ used in this paper are denoted by $He_n$ in \cite{Abramowitz:92::}.

Now, we investigate the differential Galois group $\widehat G$ of
equation~\eqref{eq:p1} over our base field   $\C(q,p)$. More precisely we
determine its dimension.  Let us recall that according to the Kolchin theorem
the dimension of the differential Galois group of an equation   is equal to  the
transcendence degree of the Picard-Vessiot extension of this field solving the
equation.  Knowing the dimension we can decide  if the group is virtually
Abelian.
\begin{lemma}
\label{l}
 If $\lambda \not\in \Z$, then $\widehat G$ is not virtually Abelian.
If $\lambda \in \Z$, then $\widehat G$ is virtually Abelian.
\end{lemma}
\begin{proof}
Let $K$ be the Picard-Vessiot extension of $\C(p)$ solving the linear system
\begin{equation}
\left\{
 \begin{split}
  & x''  +p x' + \lambda x = 0, \\
 &  u' + p u = 0,
 \end{split}\right.
\end{equation}
so
\begin{equation*}
 K := \C(p, u, x_1,x_2,x_1', x_2'),
\end{equation*}
for a basis of solutions.  Let $K_1$ be the subextension of $K/\C(p)$ generated
by $u$, i.e., $K_1 = \C(p, u)$,  and $K_2$ be the differential subextension of
$K/\C(p)$ generated by $x_1$ and $x_2$,  i.e.,
\begin{equation*}
 K_2 := \C(p,  x_1,x_2,x_1', x_2').
\end{equation*}
Note, that according to our notation $K_1=\C(p, u)=\C(p,\exp(-p^2/2))\simeq \C(q,p)$.

 We have two towers of extensions
\begin{equation}
\label{t1}
 \C(p)\subset K_1 \subset K,
\end{equation}
and
\begin{equation}
\label{t2}
 \C(p)\subset K_2 \subset K.
\end{equation}
Our aim is to determine the transcendence degree $\trdeg(K/K_1)$ of extension
$K/K_1$. This number is the dimension of the differential Galois group
$\widehat G$.

Using basic properties of the transcendence degree, see,  e.g., Chapter~8 in \cite{Lang:02::}, from the first
tower~\eqref{t1} we obtain
 \begin{equation}
  \trdeg(K/\C(p)) = \trdeg(K/K_1) + \trdeg(K_1/\C(p)).
 \end{equation}
We know that $\trdeg(K_1/\C(p))=1$ because $u=\exp(-p^2/2)$ is transcendent over
$\C(p)$.  Thus, we need to determine $ \trdeg(K/\C(p))$. But from the tower
\eqref{t2} we have
\begin{equation*}
 \trdeg(K/\C(p)) = \trdeg(K/K_2) + \trdeg(K_2/\C(p)),
\end{equation*}
and this  gives us
\begin{equation}
 \trdeg(K/K_1) = \trdeg(K/K_2) + \trdeg(K_2/\C(p)) - 1.
\end{equation}
Two cases have to be distinguished.

\begin{itemize}
\item If $\lambda$ is not integer, then by Proposition~\ref{pro:rd}, the
      transcendence degree of $K_2 / \C(p)$ is 4. So transcendence degree of
      $K/K_1$ must be greater or equal to $3$. The Galois group of this extension
      is a dimension 3 subgroup of $\mathrm{GL}(2,\C)$, so it is not virtually
      Abelian.
 \item If $\lambda$ is integer, then by Proposition~\ref{pro:s},
      $q \in K_2$. Thus, the transcendence degree of $K/K_2$ is $0$. As, by
      Proposition~\ref{pro:rd}, $\trdeg(K_2/\C(p))=2$, the transcendence degree
      of $K/K_1$ is 1. The Galois group of the variational equation is virtually
      Abelian.
\end{itemize}
\end{proof}

\begin{corollary}
\label{cor:int}
 If an eigenvalue of $V''(\vd)$ at a Darboux point is not integer, then the
 Hamiltonian system~\eqref{heq} is not integrable.
\end{corollary}

\section{Rank 4 subsystems}
\label{sec:r4}
In this section we investigate the differential Galois group over $\C(q,p)$ of
system \eqref{s12}-\eqref{s13} under assumption that $\lambda\in\Z$. At first we
show the following.
\begin{proposition}
\label{pro:s12}
If $\lambda \in \Z$, then equation \eqref{s12} has a solution in $\C(p,q)$,  and
its differential  Galois group over $\C(p,q)$ is the additive group
$G_{\mathrm{a}}$.
\end{proposition}
\begin{proof}
If we take $p$ as  an independent variable in equation~\eqref{s12}, then it becomes
equation~\eqref{eq:p1}. Now, by Proposition~\ref{pro:s} its solution $x_\lambda
\in\C(q,p)$.
We know from the proof of Lemma~\ref{l} that the dimension of the Galois group
is $1$, and this proves the second statement.
\end{proof}
\begin{proposition}
\label{pro:ga}
If $\lambda \in \Z$, and the differential Galois group over $\C(q,p)$ of
system~\eqref{s12}-\eqref{s13} is virtually Abelian, then there exists
$a\in\C$ such that the integral
\begin{equation}
\label{eq:R}
 R:=\int \left(a\frac{1}{x^2_\lambda}  - \frac{1}{q^2}x^2_\lambda \right) \rmd t=
 a \int  \frac{q}{x^2_\lambda} \rmd p - \int\frac{x^2_\lambda}{q} \rmd p,
\end{equation}
belongs to  $\C(q,p)$.
\end{proposition}
\begin{proof}
By Proposition~\ref{pro:s12}, the differential Galois group of
equation~\eqref{s12} is the additive group $G_{\mathrm{a}}$. So, we can apply
point 2 of Theorem~2.3 from \cite{Maciejewski:09::}. The condition $(\alpha)$ in
this theorem is, in our case, $R\in\C(q,p)$, for a certain $a\in\C$.
\end{proof}

\begin{lemma}
For any $a\in \C$ the integral $R(q,p)$  does not belong to the field $\C(p,q)$.
\end{lemma}
\begin{proof}
We prove our Lemma by a contradiction. Let us assume that there exists $a\in\C$
such that $R(q,p)\in\C(p,q)$.

At first we consider the case $\lambda \geq 1$. In this case the solution of
equation \eqref{s12} in $\C(p,q)$ is $x_\lambda = q H_{\lambda-1}(p)$, see
Proposition~\ref{pro:s} and thus
\begin{equation}
\label{eq:l>}
 R = \int \left( a \frac{1}{q}\frac{1}{H_{\lambda-1}(p)^2} -
 q H_{\lambda-1}(p)^2\right)\rmd p.
\end{equation}
 As a rational function of the transcendent $q=\exp(\varepsilon-p^2/2)$ with
coefficients in $\C(p)$, $R$ can be written in the form
\begin{equation*}
  R = \sum_{l=0}^m \alpha_l q^l + \frac{N(p,q)}{D(p,q)},
\end{equation*}
where $\alpha_k \in \C(p)$, $N$ and $D$ are elements of $\C(p)[q]$, and
$\deg_q(N) < \deg_q(D)$.
Moreover, this decomposition is unique.
Differentiating both sides of~\eqref{eq:R} we obtain
\begin{equation}
\label{eq:lr}
a \frac{1}{q}\frac{1}{H_{\lambda-1}^2} - q H_{\lambda-1}^2
  = \sum_{l=0}^m\left(\alpha_l' - lp\alpha_l\right) q^l + \left( \frac{N(p,q)}{D(p,q)}\right)'.
\end{equation}
As
\begin{equation*}
 0>\ord_q \frac{N}{D}:=\deg_q N - \deg_qD\geq \ord_q \left( \frac{N}{D}\right)',
\end{equation*}
we must have $m=1$, because the decomposition in both sides of~\eqref{eq:lr} is
unique. In particular, we have
\begin{equation}
\label{a1}
 \alpha_1' - p\alpha_1 = -H_{\lambda-1}^2 .
\end{equation}
Because this equation is regular, $\alpha_1 \in \C[p]$, and it can be written as a
linear combination of  Hermite polynomials
\begin{equation*}
 \alpha_1 = \sum_{n=0}^N \gamma_n H_n, \mtext{where} \gamma_n\in\C.
\end{equation*}
Hermite polynomials satisfy the following relation
\begin{equation*}
 H_{n+1}(p)=pH_n(p)-H_n'(p),
\end{equation*}
so we have
\begin{equation*}
  \alpha_1' - p\alpha_1= \sum_{n=0}^N \gamma_n(H'_n-pH_n)= -\sum_{n=0}^N \gamma_n H_{n+1}.
\end{equation*}
Thus, we can rewrite~\eqref{a1} in the form
\begin{equation*}
 H_{\lambda-1}^2  = \sum_{n=0}^N\gamma_n H_{n+1},
\end{equation*}
and so
\begin{equation*}
 \int_{-\infty}^\infty \rme^{-p^2/2} H_{\lambda-1}^2(p)\,\rmd p =
 \sum_{n=0}^N\gamma_n\int_{-\infty}^\infty \rme^{-p^2/2} H_{n+1} (p)\,\rmd p.
\end{equation*}
But this gives a contradiction
\begin{equation*}
 \sqrt{2\pi} (\lambda-1)!=0,
\end{equation*}
because
\begin{equation*}
 \int_{-\infty}^\infty \rme^{-p^2/2} H_{n+1} (p)\,\rmd p=
 \int_{-\infty}^\infty \rme^{-p^2/2} H_{n+1} (p)H_0(p)\,\rmd p=0,
\end{equation*}
for $n\geq 0$. For $\lambda\geq 1$ our lemma is proved.

In the case $\lambda\leq 0$  the solution of equation  \eqref{s12} in
$\C(p,q)$ is $x_{\lambda} = H_{\lambda}(-ip)$, so integral~\eqref{eq:R} reads
\begin{equation}
 R=\int\left( \frac{aq}{H_{-\lambda}(-\rmi p)^2}-
 \frac{1}{q}H_{-\lambda}(-\rmi p)^2\right)\rmd p.
\end{equation}
If we set
\begin{equation}
 v=\rmi p, \quad u=\rme^{\varepsilon-v^2/2}=\rme^{2\varepsilon}/q, \quad \widetilde \lambda=  1-\lambda, \mtext{and} \widetilde a =a \rme^{4\varepsilon} ,
\end{equation}
then we transform the considered integral into the following one
\begin{equation}
 R = \rmi \rme^{-2\varepsilon}\int \left( \frac{\widetilde a }{uH_{\widetilde\lambda-1}(v)^2} -  uH_{\widetilde\lambda-1}(v)^2\right)\rmd v.
\end{equation}
But this integral, after renaming variables, is proportional to that one already
considered for $\lambda\geq1$, see~\eqref{eq:l>}. So, it is not rational and
this finishes the proof.
\end{proof}
As  corollaries we have.
\begin{corollary}
 If $\lambda\in\Z$, then the differential Galois group over $\C(q,p)$ of
 system~\eqref{s12}-\eqref{s13} is not virtually Abelian.
\end{corollary}
\begin{corollary}
\label{cor:block}
If $V''(\vd)$ at a Darboux point has a Jordan block with integer eigenvalue,
 then Hamiltonian system~\eqref{heq} is not integrable.
\end{corollary}
Now, our main result given by Theorem~\ref{thm:we} follows directly from
Corollary~\ref{cor:int} and Corollary~\ref{cor:block}.
\section{Examples}
\label{sec:ex}
We consider in details case $n=2$. Our aim is to give a full characterisation of
homogeneous potentials $V\in\C(q_1,q_2)$ of degree $k=0$, which satisfy the
necessary conditions of Theorem~\ref{thm:we}.

Darboux points of $V$ are non-zero solutions of equations
\begin{equation}
  \label{eq:aux}
    \dfrac{\partial V}{\partial q_1}=q_1\qquad
 \dfrac{\partial V}{\partial q_2}=q_2.
\end{equation}
As it was explained in \cite{Maciejewski:04::e} it is convenient to consider
Darboux points as points in the projective line $\CP^1$. Let $z=q_2/q_1$,
$q_1\neq 0$, be the affine coordinate on $\CP^1$. Then, we can rewrite
system~\eqref{eq:aux} in the form
\begin{equation}
 v'(z)z=-q_1^2 , \qquad  v'(z) = zq_1^2,
\end{equation}
where $v(z):=V(1,z)$.  From the above formulae it follows that $z_\star$ is a Darboux point of $V$, if and only if  $z_\star \in\{-\rmi, \rmi\}$, and $v'(z_\star)\neq 0$. Thus, the location of Darboux points does not depend on the form of potential!

If $z_\star$ is the affine coordinate of a Darboux point $\vd$ of $V$, then the Hessian matrix $V''(\vd)$ expressed in this coordinate  has the form
\begin{equation}
 V''(\vd)=\begin{bmatrix}
             -v''(z_\star) x_\star^{-2} - 2 & -[v'(z_\star) +z_\star v''(z_\star)]x_\star^{-2} \\
 -[v'(z_\star) +z_\star v''(z_\star)]x_\star^{-2}  &  v''(z_\star)x_\star^{-2}
          \end{bmatrix}
\end{equation}
where
\begin{equation*}
 x_\star^2 = -  v'(z_\star)z_\star = v'(z_\star)/z_\star.
\end{equation*}
Vector $\vd$ is an eigenvector of $V''(\vd)$ with the corresponding eigenvalue
$\lambda=-1$. As $\tr V''(\vd)=-2$, $\lambda=-1$ is the only eigenvalue of
$V''(\vd)$. Thus the first condition of Theorem~\ref{thm:we} is satisfied. If
$V''(\vd)$ is semi-simple, then it is diagonal. Hence the second condition of
Theorem~\ref{thm:we} is satisfied iff
\begin{equation}
\label{sem}
 v'(z_\star) +z_\star v''(z_\star)=0.
\end{equation}
Let us apply the above criterion for potential
\begin{equation*}
 V =\frac{ q_2}{q_1^3}(q_2-a q_1)(q_2-b q_1) \mtext{where} a\neq b,
\end{equation*}
assuming that it has two Darboux points with affine coordinates $\pm\rmi$. An
easy calculations shows that condition~\eqref{sem} is satisfied for
$z_\star=\pm\rmi$, iff
\begin{equation*}
 V = \frac{q_2}{q_1^3}(9q_1^2+ q_2^2).
\end{equation*}
We perform a  direct search of the first integrals  for this potential. We looked for first integrals   which are polynomial with respect to the momenta,  and we assumed that he coefficient of these polynomials  are differentiable functions.  In this way we checked that the potential does not admit any additional first integral which is a polynomial with respect to the momenta of degree less than five. So it is not clear if the potential is integrable or not. Anyway we must remember that Theorem~\ref{thm:we} gives only necessary, not sufficient conditions for the integrability.
\section*{Acknowledgements}

This research has been partially supported  by grant No. N N202 2126 33 of Ministry of Science and Higher Education of Poland,  by  UMK grant  414-A, and by EU funding for the Marie-Curie Research Training Network AstroNet.

\newcommand{\noopsort}[1]{}\def\cprime{$'$} \def\cprime{$'$} \def\cprime{$'$}
  \def\cydot{\leavevmode\raise.4ex\hbox{.}} \def\cprime{$'$}

\end{document}